\documentclass[11pt]{article}
\usepackage{natbib}
\usepackage[american]{babel}
\usepackage{amsmath}
\usepackage{amssymb}
\usepackage{mathrsfs}
\usepackage{amsbsy}
\usepackage{subfig}
\usepackage{epsfig}
\usepackage{boxedminipage}
\usepackage{graphicx}
\usepackage{lineno}
\usepackage{setspace}

\newcommand{\be}{\begin{equation}}
\newcommand{\ee}{\end{equation}}

\begin{document}

\title{Error Control of Iterative Linear Solvers for Integrated Groundwater Models}

\author{
Matthew F. Dixon\thanks{Department of Computer Science, One Shield's Avenue, 
University of California, Davis, CA 95616}~\thanks{Corresponding 
author Email: mfdixon@ucdavis.edu, Tel: 530-754-9586, Fax: 530-752-4767.} ~
Zhaojun Bai\footnotemark[1] ~
Charles F. Brush\thanks{Bay-Delta Office, Department of Water 
Resources, 1416~ 9th Street, Sacramento, CA 95814} ~ 
Francis I. Chung\footnotemark[3] ~ \\ 
Emin C. Dogrul\footnotemark[3] ~
Tariq N. Kadir\footnotemark[3]
}




\maketitle


\abstract{An open problem that arises when using modern iterative linear solvers, such as the preconditioned conjugate gradient (PCG) method or Generalized Minimum RESidual method (GMRES) is how to choose the residual tolerance in the linear solver to be consistent with the tolerance on the solution error.  This problem is especially acute for integrated groundwater models which are implicitly coupled to another model, such as surface water models, and resolve both multiple scales of flow and temporal interaction terms, giving rise to linear systems with variable scaling.  

This article uses the theory of 'forward error bound estimation' to show how rescaling the linear system affects the correspondence between the residual error in the preconditioned linear system and the solution error. Using examples of linear systems from models developed using the USGS GSFLOW package and the California State Department of Water Resources' Integrated Water Flow Model (IWFM),  we observe that this error bound guides the choice of a practical measure for controlling the error in rescaled linear systems. It is found that forward error can be controlled in preconditioned GMRES by rescaling the linear system and normalizing the stopping tolerance. We implemented a preconditioned GMRES algorithm and benchmarked it against the Successive-Over-Relaxation (SOR) method.  Improved error control reduces redundant iterations in the GMRES algorithm and results in overall simulation speedups as large as 7.7x. This research is expected to broadly impact groundwater modelers through the demonstration of a practical approach for setting the residual tolerance in line with the solution error tolerance.



}


\section*{Introduction}

As the groundwater model infrastructure grows to comprehensively and accurately resolve hydrological processes so too does the need for improved solver technology to address new modeling features and take advantage of faster computers. Consequently, there are now a wide range of iterative linear solvers available in groundwater modeling packages. Examples include the preconditioned conjugate gradient (PCG) method \citep{Hill90}, the link-algebraic multi-grid (LMG) Package \citep{Mehl01}, the algebraic multi-grid (AMG) solver and the generalized conjugate gradient method. The first three of these are provided with MODFLOW-2005 \citep{Harbaugh05} and the latter is provided in SEAWAT \citep{Guo02}, which couples MODFLOW-2005 and MT3DMS \citep{Zeng99} to simulate ground water flow with variable density and temperature. 


Iterative linear solvers can be broadly categorized into modern, projection based, solvers or classical (stationary) solvers. Both can be further categorized into solvers for symmetric or non-symmetric linear systems. The PCG method remains one of the most competitive modern solvers widely used for groundwater modeling. Whilst this solver is designed for symmetric matrices only, it shares a common property with numerous other modern solvers, namely that the most robust stopping criteria is based on the residual error\footnote{The stopping criteria in the Generalized Minimum RESidual method (GMRES) can only be based on an estimate of the residual error norm.}. It is, however, the solution error which is the most physically relevant error measure - solving a linear system to within 1\% accuracy in the groundwater heads is more relevant than to within 1\% accuracy of the linear system residual. An open problem is how to control the solution error given a stopping tolerance on the residual error, which behaves very differently when the linear system is ill-conditioned \citep{axel01}.

The issue of choosing a residual based stopping criteria becomes even more apparent when a projection based iterative solver is embedded in an iterative linearization method (such as the Picard or Newton method), which is generally necessary for modeling saturated groundwater flow. A relative error based stopping tolerance for the linearization procedure is generally not compatible with a residual based stopping criteria in the linear solver. In this situation, the only remedy is a judicious choice of stopping tolerance on the residual error, often based on trial and error, which ensures the linearization procedures quickly converges without redundant linear solver iterations. The problem with this approach is that to achieve a target solution error throughout the simulation, the corresponding solver error tolerance may need to change as the linear systems are forced with temporal effects such as pumping, stream seepage or wetting and drying of the aquifers. This makes the task of choosing the stopping tolerance even more difficult and there may be little choice but to over-specify the tolerance, before running the simulation, at the expense of excessive linear iterations and performance lagging.

The multiple scales of flow in groundwater models coupled with surface water or contaminant transport such as GSFLOW \citep{Markstrom08} or MT3DMS further highlight the issues associated with linear solver error control. \citet{Blom93} consider the scaling issues arising between the solution components of a model for brine transport in groundwater flow. They use a weighted norm in the linear solver stopping criteria in order to ensure that each solution component is solved to its corresponding data accuracy. \citet{Blom93} further considered the influence of the linear solver error on the convergence of a Newton-type method. They propose a fixed bound on the linear solver error, referred to as the \emph{forward error} which is shown to be inversely proportional to the maximum number of Newton iterations. They don't, however, show how this forward error is related to the stopping tolerance and thus it is unclear as to how this approach can be efficiently implemented in practice.

\emph{This paper shows how re-scaling poorly scaled linear systems is a key step towards improving control of the forward error, since the (normalized) stopping tolerance on the residual norm becomes a practical proxy for the upper bound on the forward error. Improved control eliminates unnecessary linear solver iterations without compromising the desired accuracy of the solver. For saturated flows, too much error can reduce the convergence rate of the linearization procedure and even prevent convergence.} 

We demonstrate the practical benefits of rescaling the linear system in two different integrated surface water and groundwater models - the GSFLOW package and the Integrated Water Flow Model (IWFM) \citep{IWFM07}. IWFM is a water resources management and planning tool which simulates groundwater, surface water and stream-groundwater interaction. This model is currently being used by the State of California Department of Water Resources in computationally demanding long-time high resolution applications such as assessing the impact of climate change on water resources and the analysis of different conjunctive use scenarios across California. IWFM uses an implicitly formulated Galerkin finite element method over a non-uniform areal 2D mesh to simulate the nonlinear groundwater head dynamics in multi-layer aquifers. 

Whilst both GSFLOW and IWFM implicitly couple the surface water flow with the groundwater flow, IWFM combines both sets of flow equations into an integrated linear system with a non-symmetric matrix and hence the PCG solver can no longer be used. The need for non-symmetric matrix solvers is a growing trend in groundwater modeling. As previously mentioned, SEAWAT uses a generalized conjugate gradient method which is suitable when the matrices are nearly symmetric.  The need for faster local converge rates than attainable using Picard methods motivates the use of Newton-type methods for saturated ground water models. By using full upstream weighting of the saturated thicknesses to compute the inter-cell conductances in MODFLOW-2005, Newton-type methods give non-symmetric linear systems and require different solvers and settings. \citet{Mehl06} concludes that further exploration into the different solvers and their settings is needed before this approach can widely catch on.  To partially address this call for further exploration, we provide an overview of a promising modern iterative solver for a wider class of non-symmetric matrices, referred to as the Generalized Minimum RESidual (GMRES) method \citep{Saad86}, which is known to be superior to the classical Successive Over Relaxation (SOR) method \citep{Vorst00,Pad08}.

\subsection*{Overview}
The next Section describes the properties of linear systems that are most performance critical for any iterative solver and identifies the prevailing features of numerous case studies to explain why integrated surface water and groundwater models are more difficult to solve. Almost all of the linear systems that we consider have non-symmetric coefficient matrices and exhibit scaling difficulties. We then briefly review the mathematical formulation of GMRES, widely used for scientific computations involving non-symmetric linear systems. This solver has only received marginal coverage within the groundwater modeling literature and we point out some of its salient features and practical implications for groundwater modeling. We also provide the algorithm parameter settings that lead to best performance and benchmark it against the classical SOR method. The key contribution of this paper is to provide insight into how to choose the residual stopping tolerance in any modern iterative linear solver (which uses residual based stopping criteria), not just the GMRES method. Whilst this is an open problem, we show that error control in GMRES can be improved by rescaling the linear system if it is poorly scaled, as may arise when groundwater models are integrated with surface water flow models such as IWFM and GSFLOW. 

\section*{Profile of the Linear System} \label{sect:pro}
At each time step in a saturated groundwater model simulation, a linearization procedure such as the Picard or Newton iterative methods (see for example \citep{Mehl06}) solves 
the system of saturated groundwater flow equations 
$$F(H^{k+1})=0$$ 
in which $H^{k+1}$ is the vector of unknown multiple-layer aquifer heads over a 2D bounded domain at iteration $k+1$. This definition is general enough to include implicitly coupled integrated groundwater and surface models (such as IWFM), which also include stream and lake surface elevations in the vector of unknowns. For ease of exposition we present the linear system in canonical form by denoting the difference vector $x=H^{k+1}-H^{k}$ without an iteration index, the Jacobi (or approximate Jacobian) matrix $A = \nabla F(H^k)$ with 
elements $a_{ij}=\frac{\partial F_i}{\partial H^k_j}$ and the right-hand side vector $b=-F(H^k)$ to give
\be \label{eqn:linsys}
Ax=b, \qquad A\in\mathbf{R}^{N\times N},~x,b\in\mathbf{R}^N,
\ee
where $A$ is a positive definite\footnote{A matrix $A$ is positive definite if $x^TAx\geq 0$ for all real $x\neq 0$.} square matrix. This paper will just consider the efficient iterative approximation of linear system \ref{eqn:linsys} for the case when $A$ is non-symmetric. However, the issue of accuracy control of iterative solvers discussed later in this paper also applies to any positive definite matrix (symmetric or not). 

Table \ref{tab:param} shows the performance critical properties 
of the coefficient matrix $A$ for seven 
datasets arising in various groundwater packages. The first four data sets are from applications using IWFM. HCMP is a synthetic hydrological comparison dataset, C2VSIM and C2VSIM9 are from respective three and nine aquifer layer integrated models of California's Central Valley and BUTTE is from a high resolution integrated model of BUTTE county. INCLINE and SAGEHEN are from GSFLOW models of the Incline and Sagehen water creeks respectively (see \citet{Markstrom08} for details of the Sagehen water creek example). NAC is from a two-layer model of Nacatoch Aquifer \citep{texas09} which uses MODFLOW-2000. The INCLINE and the NAC datasets were produced using a Newton Method in a development version of the MODFLOW package. 

Dimension $N$ is the size of the matrix and NNZ denotes 
the number of non-zero elements. 
Each matrix is sparse and lacks any block structure. 
Sparsity is the percentage of the elements in a matrix which are non-zero.
Normality is the relative measure $\|AA^*-A^*A\|/\|A\|^2$ 
which is zero when $A$ is symmetric. $\kappa(A)$ is the estimated 
condition number\footnote{$\kappa(A)$ is estimated using 
the {SuperLU} \citep{SuperLU} routines {\tt dgscon} and {\tt dlangs}.} 
of $A$ and is a measure of sensitivity of the linear system and 
the convergence rate of iterative solvers. 

\begin{table}
\begin{center}
\begin{tabular}{|l|l|l|l|l|l|l|l|} \hline
 & \multicolumn{4}{|c|}{IWFM} & \multicolumn{3}{|c|}{GSFLOW / MODFLOW}\\
 \hline
             & HCMP &  C2VSIM  & C2VSIM9 & BUTTE & INCLINE & NAC & SAGEHEN\\ \hline
  Dimension  & 46460  & 4630  &     12988  & 34683 & 42820  & 75319  &  6784 \\
  NNZ      & 479246  &   41616  &     125616 & 188006 &  261696 & 433791& 31504 \\
  Sparsity(\%)  &  0.0220   & 0.194  &    0.0744 & 0.0156 & 0.0143 & 0.00765 & 0.0685\\
  Normality  &   0.271  & 0.222  &     0.908 &  0.199 & 0.00730 &  0.0114 & 0 \\
  $\kappa(A)$  & 3.09E6   &  2.54E11  & 5.13E6 & 1.95E9& 1.50E8 & 7.25E4 & 6.71E8\\ 
  \hline
\end{tabular}
\caption{Linear solver performance critical properties for seven different datasets taken from applications using IWFM, GSFLOW and MODFLOW.}  
 \label{tab:param}
 \end{center} 
 \end{table}

Figure 1 illustrates scaling issues arising in the coefficient matrices from integrated ground water and surface water models. Figures 1a and 1b show the sparsity patterns of the C2VSIM and SAGEHEN coefficient matrices respectively. Each Figure has been separated into distinctive zones for illustrative purposes and a color scheme arbitrarily differentiates scale. In Figure 1a, the sparsity structure in the uniform $3\times 3$ grid corresponds to the three aquifer layers of the C2VSIM model and their interactions with each other. For example, the middle layer interacts with the layer above and below and thus exhibits a diagonal band for each together with a central band for the convection and diffusion within the layer. The bottom and top layers only interact with one other aquifer layer and thus only exhibit two bands. The non-uniform micro-structure within the central band typifies that of an unstructured Galerkin finite element groundwater water flow model. The upper left-hand zone of Figure 1a corresponds to the stream nodes and the zones to the right and below correspond to the stream-groundwater interaction terms with the top level aquifer. The $2\times2$ uniform grid in Figure 1b corresponds to two aquifer layers of the SAGEHEN water creek model and their interactions with each other. Both Figures illustrate the scaling issues. In Figure 1a, matrix elements whose absolute values are above and below an arbitrary threshold of $\mathcal{O}(10^6)$ are shown in red and blue, respectively, whilst in Figure 1b, the threshold is $\mathcal{O}(10)$. 


Figures 1c and 1d respectively show the corresponding graphs of the C2VSIM and SAGEHEN coefficient matrix element sizes to further illustrate the scaling issues. Figure 1c is split into the rows corresponding to the surface water and top aquifer layer by the vertical red line. The vertical axis is a power scale for the absolute matrix element sizes in each of the first one thousand rows whose indices are shown on the horizontal scale. The non-zero matrix elements corresponding to the stream nodes in the C2VSIM model are not only much sparser than those corresponding to the aquifer nodes, but exhibit a broader range of absolute values. Figure 1d is split into the rows corresponding to the top and bottom aquifer layers of the SAGEHEN model by the vertical red line.

\begin{figure}[!ht]
\vskip.3cm
\epsfxsize=10cm
\centering
\subfloat[][]{\includegraphics[angle=0.0,width=0.5\textwidth, height=0.3\textheight]{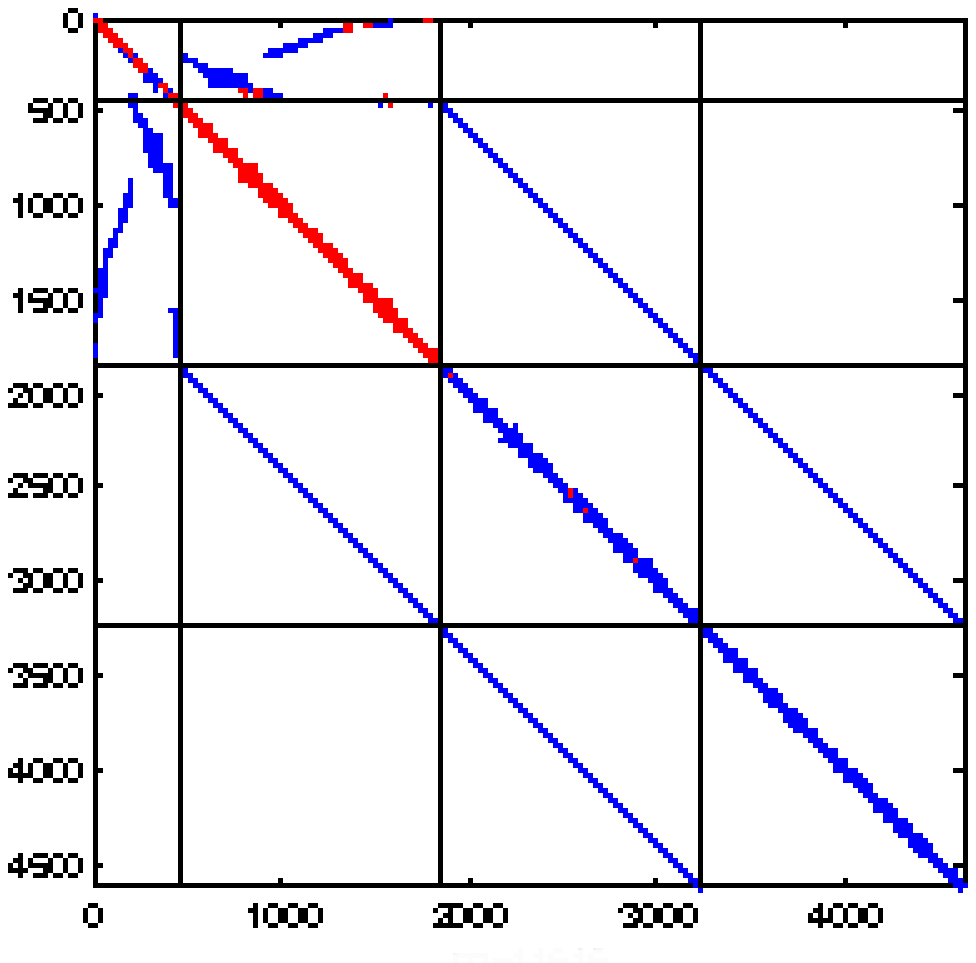}}
\subfloat[][]{\includegraphics[angle=0.0,width=0.5\textwidth, height=0.3\textheight]{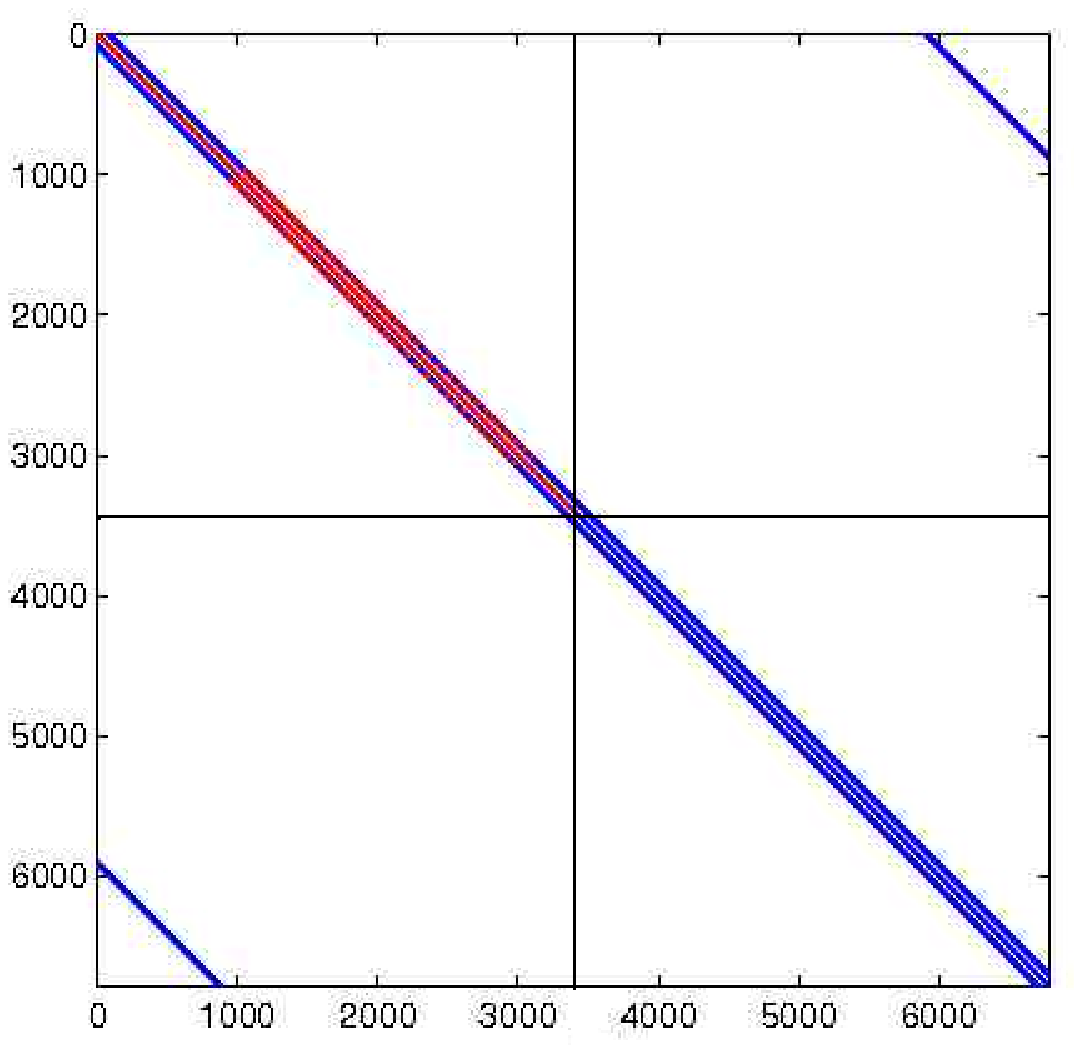}}\\

\subfloat[][]{\includegraphics[angle=0.0,width=0.4\textwidth, height=0.3\textheight]{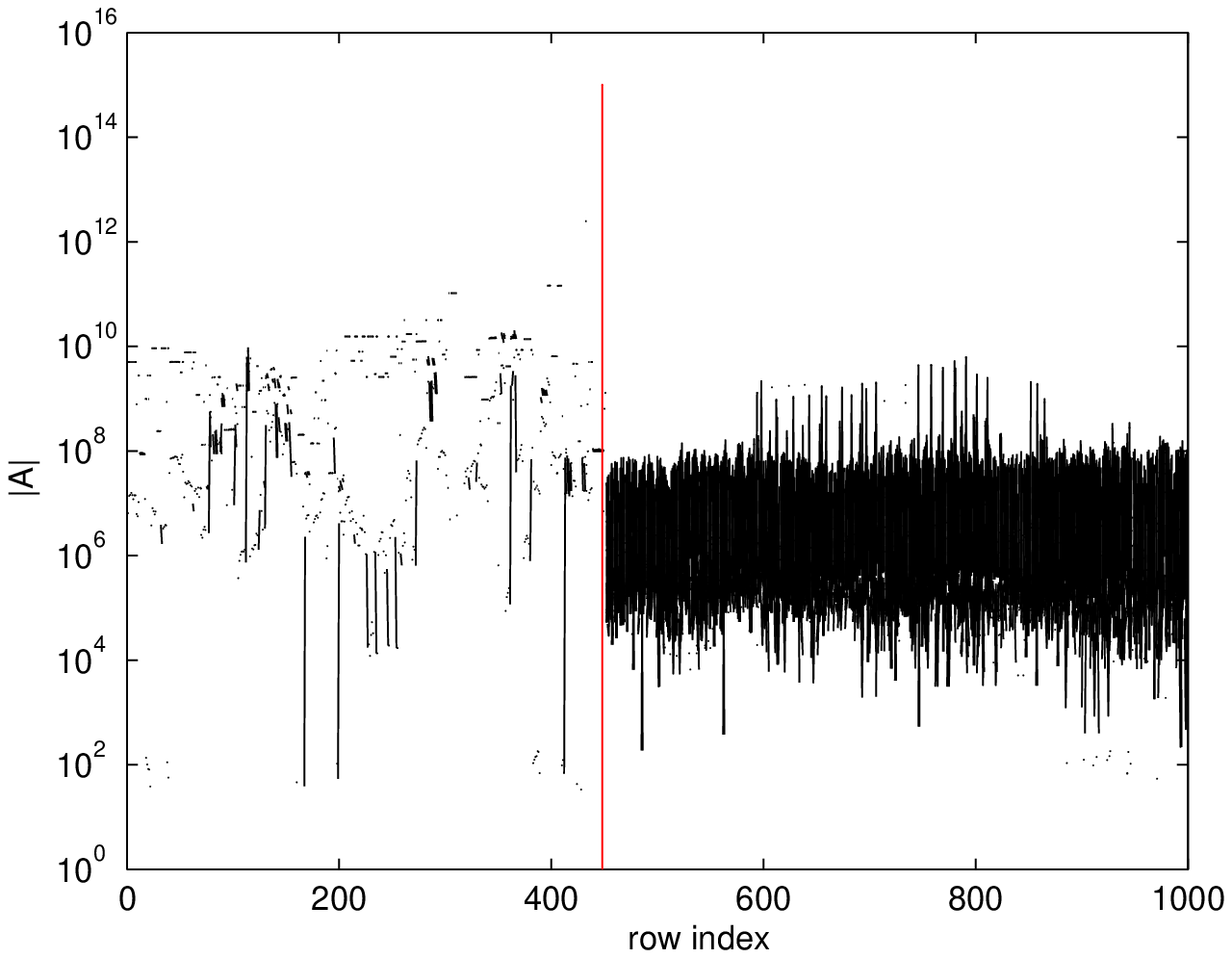}}
\qquad \qquad
\subfloat[][]{\includegraphics[angle=0.0,width=0.4\textwidth, height=0.3\textheight]{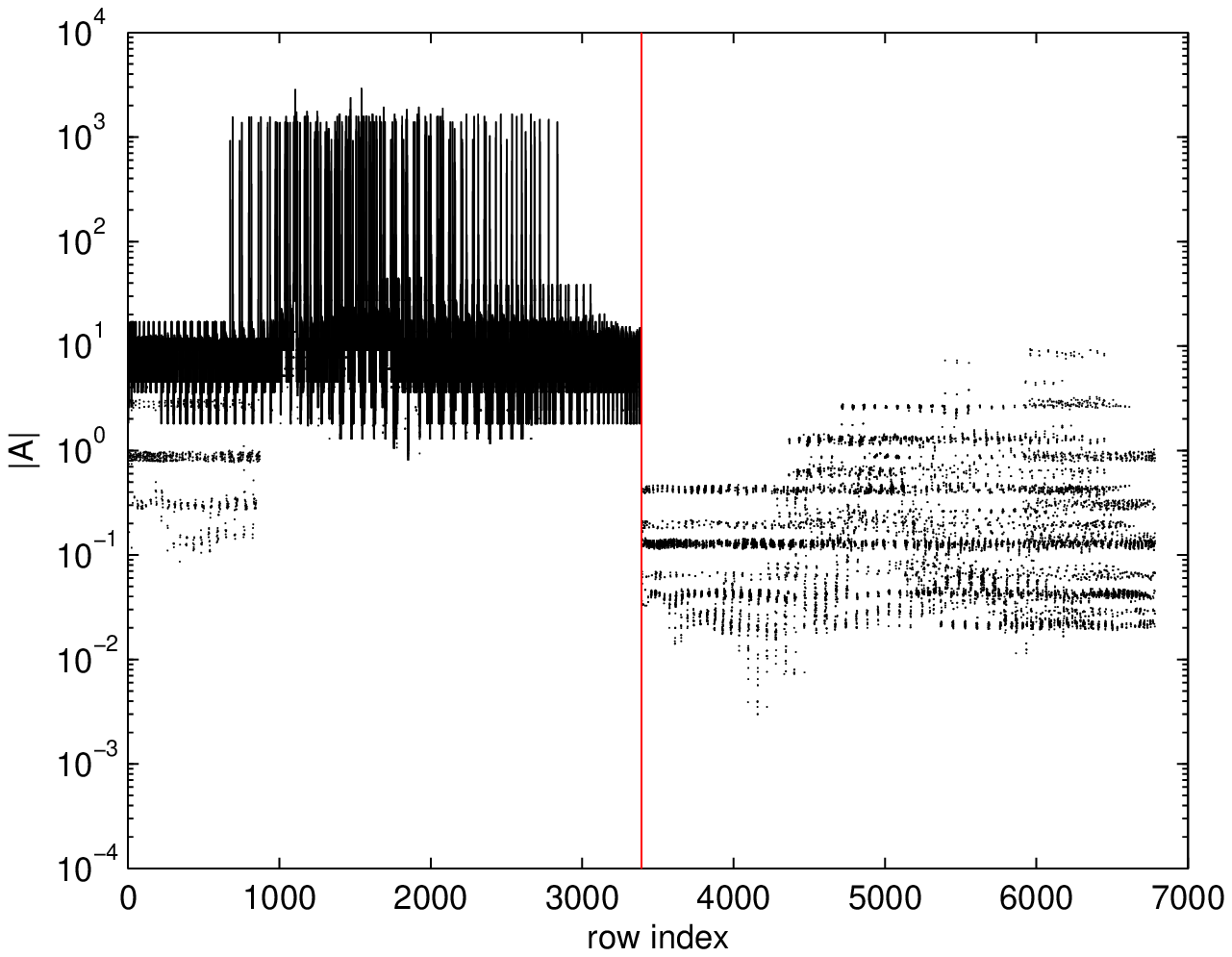}}\\
\caption{The sparsity pattern of the (a) C2VSIM coefficient matrix and (b) SAGEHEN coefficient matrix. (c) A log scale of the absolute value of the elements in the (c) C2VSIM coefficient matrix against the row index (which has been truncated to highlight the interface region) and the (d) SAGEHEN coefficient matrix against the row index.}
\end{figure}

\bigskip 


\paragraph{Objectives} Having illustrated some of the scaling issues arising in the coefficient matrices from IWFM and GSFLOW, our objectives are to 

\begin{itemize}
\item Solve the linear system \eqref{eqn:linsys} iteratively to a desired level of accuracy in the solution; 
\item Assess the impact of the scaling issues on solver robustness and performance; and
\item Provide practical techniques and optimal parameterizations that readers can either implement or use to improve solver robustness and performance.
\end{itemize}

The following Section describes the use of GMRES with preconditioning to efficiently solve  \eqref{eqn:linsys}. This algorithm terminates when an estimate of the preconditioned linear system residual norm is below a tolerance $\tau$. \emph{The problem is how to choose $\tau$ so that the solution is to within the desired accuracy}. If the estimated condition number of the preconditioned coefficient matrix is known, then this task becomes trivial. However, the coefficient matrices are poorly conditioned and may vary considerably throughout the simulation. The matrices are also too large for efficient estimation of the condition number prior to solving each preconditioned linear system.  

\emph{The solution is to rescale the linear system so that the bound of the relative forward error of the linear
system is about the same order of magnitude as the stopping tolerance normalized by the right-hand-side vector}. This avoids the computationally prohibitive step of having to estimate the condition number of the preconditioned coefficient matrices before each solve and instead provides a simple approach with negligible additional computation.



\bigskip 
 
\section*{The GMRES Algorithm and Preconditioning} \label{sect:gmres}
The Generalized Minimum RESidual (GMRES) method is a Krylov subspace projection method for solving the linear system \eqref{eqn:linsys} based on taking the pair of projection subspaces
\be
\mathcal{W}=\mathcal{K}_m(A,r_0) \qquad\textrm{and}\qquad \mathcal{V}=A\mathcal{W},
\ee
where $\mathcal{K}_m(A,r_0)$ is a Krylov subspace defined as
\be
\mathcal{K}_m(A,r_0):=\text{span}\{r_0,Ar_0,A^2 r_0,\dots,A^{m-1}r_0\} 
\ee
and $r_0=b-Ax_0$ with an initial approximate solution $x_0$. 
An approximation solution 
$\hat{x}\in x_0 + \mathcal{W}$ has the form $\hat{x}=q_{m-1}(A)r_0$ 
and $A\hat{x} - b \perp \mathcal{V}$, 
where $q_m(A)$ is a matrix polynomial of degree $m$. 

GMRES first uses an Arnoldi procedure to build an orthonormal matrix 
$V_m=[v_1,v_2,\dots,v_m]$ whose column vectors span the subspace 
$\mathcal{W} = \mathcal{K}_m(A, r_0)$. 
In matrix notation, 
the Arnoldi procedure can be expressed by the following governing 
equation 
\be
AV_m=V_{m+1}\widehat{H}_m,
\ee
where $\widehat{H}_m:=[H^T_m, h_{m+1,m}e_m]^T$, and 
$H_m$ is $m \times m$ upper Hessenberg matrix, 
and $e_m$ is the last column $m$-vector 
in the identity matrix $I_m$. 
An iterative solution to the linear system \eqref{eqn:linsys} can be 
written in the form $x_m=x_0 + V_my_m$, where the $m$-vector 
$y_m$ is the solution to the least squares problem 
\be
y_m = \arg \min_{y} ||r_m|| = \arg \min_{y} || \beta e_1 - \hat{H}_m y||,
\ee
which minimizes the residual. 
Thus GMRES finds the best $x_m$ which minimizes the residual $r_m$ by 
reducing $A$ to $\widehat{H}_m$ using the orthonormal bases 
$V_m$ and $V_{m+1}$. 
We refer the reader to \citet{Demmel97,Saad00} 
for a more detailed explanation of the GMRES method. 
GMRES($m$) is a memory efficient and more stable variant of GMRES, 
which resets the algorithm after $m$ iterations by setting $x_0=x_m$ 
so that the memory requirements are $\mathcal{O}(mN)$. 
$m$ is typically set to between $10$ and $20$.

Preconditioning is known as the 'determining ingredient' in the 
success of the GMRES and other iterative methods for solving
large scale problems.  The convergence rate and computational 
cost of solving the preconditioned linear system 
$$M^{-1}Ax=M^{-1}b$$ 
depend on the choice of the preconditioner $M$. 
The choice of $M$ is typically inferred from experience which tells us that the form of $M$ should (i) ensure that $\kappa(M^{-1}A)\ll\kappa(A)$ and (ii) be computationally inexpensive to solve $My=Ax$ for $y$ given a vector $Ax$.

For GMRES, an ideal choice is typically one in which $M^{-1}A$ is close 
to normal and whose eigenvalues are tightly clustered around 
some point away from the origin. 
The \emph{incomplete LU decomposition} (ILU)  is 
a popular preconditioner \citep{Saad00}. 
For example, considering an $LU$ decomposition 
$A = LU$, where $L$ is a unit lower triangular matrix and 
$U$ is an upper triangular matrix. 
Replacing non-zero elements of $L$ and $U$ outside the sparsity pattern 
of $A$ with zero elements gives the incomplete 
factors $\widehat{L}$ and $\widehat{U}$.  
An ILU preconditioner is then formed  by 
setting $M=\widehat{L}\widehat{U}$. 

A high-level description of the preconditioned GMRES($m$) method 
is provided below.

\begin{boxedminipage}[b]{\linewidth}
\begin{tabbing}
mmmmmmm\=mm\=mm\=mm\=mm\=mm\=   \kill
\>{\sc Preconditioned GMRES($m$)} \\
\>Input: $A, M, b, x_0, m, \tau $ \\ 
\>Output: $x_m, \gamma_{m+1,m}$ \\ 
\>1.\> compute $r_0 = M^{-1}(b - A x_0)$, $\beta = \|r_0\|_2$ and
       $v_1 := r_0/\beta$ \\
\>2.\> for $j = 1, 2, \ldots, m$ do\\
\>3.\>  \> solve $ Mw= A v_j $ for $w$\\
\>4.\>  \> for $i = 1, 2, \ldots, j$ do \\
\>5.\>  \>  \> $h_{ij} = v^T_i w $ \\
\>6.\>  \>  \> $w := w - h_{ij} w_i$ \\
\>7.\>  \> end do  \\
\>8.\>  \> compute $h_{j+1,j} = \|w\|_2$ and $v_{j+1} = w/h_{j+1,j}$   \\
\>9.\> end for \\
\>10.\> 
find $y_m$ so that 
$\gamma_{m+1,m}:= \|\beta e_1 - \widehat{H}_m y_m\|
= \arg \min_{y} \|\beta e_1 - \widehat{H}_m y\|$\; \\ 
\>11.\> $x_m = x_0 + V_m y_m$  \\
\>12.\> If $\gamma_{m+1,m}\leq\tau$, stop, else 
        goto Line 1 with $x_0=x_m$ 
\end{tabbing}
\end{boxedminipage}

Note that on each iteration of the PGMRES($m$) algorithm, 
the linear system $Mw=Av_k$, where $v_k\in V_k$, is 
solved for the vector $w$. When $M$ is an ILU factorization, 
$M = \widehat{L}\widehat{U}$, 
$w$ is determined by forward and back substitutions: 
\be
\widehat{L}z=Av_k \qquad \text{and} \qquad \widehat{U}w=z,
\ee 
where $\widehat{L}$ and $\widehat{U}$ are respectively a 
unit lower and upper triangular matrix with $2\cdot{\tt Lfil} +1 $ 
entries per row. The level of fill-in, {\tt Lfil}, is typically 
chosen to be between 5 and 10. The PGMRES($m$) 
algorithm terminates when the estimated residual norm 
$\gamma_{m+1,m}:=\| \beta e_1 - \hat{H}_m y_m\|$ satisfies 
the stopping criteria $\gamma_{m+1,m}\leq \tau$.

\section*{Scaling and Error Control}\label{sect:scaling}

Each stopping tolerance $\tau$ for 
the approximate solution $\widehat{x}$ computed by the
PGMRES iteration can be associated 
with a corresponding estimate of the upper bound  
on the relative forward error norm 
$\|\widehat{x}-x\|/\|x\|$ of the linear system \eqref{eqn:linsys}. 
Whilst the acceptable upper bound on the relative forward error norm 
is determined from the data accuracy, the corresponding tolerance 
can not easily be implied from the data accuracy.  
Misspecification of the tolerance can result in either 
over-solution of the linear system \eqref{eqn:linsys} 
or unacceptably high forward error with respect to the data accuracy,
especially when 
the coefficient matrix is poorly conditioned and scaled.

To reduce the difference between the estimated upper bound on 
the relative forward error norm and the residual norm, 
we introduce a diagonal scaling matrix $D$ so that the 
preconditioned linear system becomes
\be
M^{-1}D^{-1}Ax=M^{-1}D^{-1}b. 
\ee
The associated residual $\widehat{r}$ for
an approximate solution $\widehat{x}$ is defined as
$\widehat{r} =M^{-1}D^{-1}(b-A\widehat{x})$. 
By a standard perturbation analysis \citep{Demmel97,Higham02},  
the upper bound, denoted as ${\tt Ferr}$, 
on the relative forward error norm can be given 
in terms of the norm of the residual $\widehat{r}$:  
\be \label{eq:res}
\frac{\|\widehat{x} - x\|}{\|x\|} \leq 
\kappa(M^{-1}D^{-1}A)\frac{\|\widehat{r}\|}{\|M^{-1}D^{-1}b\|}  
\equiv {\tt Ferr},
\ee
where $\kappa(M^{-1}D^{-1}A^{-1})$ is the condition number
to characterize the 
difference between the relative forward error norm and the ratio of 
the residual norm to the right-hand side vector norm $\|M^{-1}D^{-1}b\|$. 
The condition number can not in general be efficiently obtained 
during simulation due to the cost and dynamic nature of 
the linear systems and thus prohibits the evaluation of ${\tt Ferr}$ 
explicitly as the stopping criteria for the PGMRES$(m)$ iteration. 
By choosing $D$ as the sum of row elements\footnote{This choice of scaling is referred to as row \emph{equilibration}. 
$e$ denotes the unit vector of length $N$.}
\be
D=\text{diag}(|Ae|_1, |Ae|_2,\dots, |Ae|_N),
\ee
we both normalize $A$ and minimize the condition number of 
$D^{-1}A$ \citep{Higham02}  thereby significantly sharpening the 
difference between the estimated upper bound on the relative 
forward error norm and the residual norm. 
\emph{Rescaling is a key step towards control of the 
forward error, since the normalized stopping tolerance $\epsilon$, where $\tau=\epsilon ||b||_*$, 
on the estimated residual norm becomes a practical proxy 
for the upper bound on the relative forward error norm \footnote{Following \citet{Blom93}, we build the scaling into the norm $||\cdot||_*:=||D^{-1}\cdot||$.}}.  

To illustrate this property, let us examine the C2VSIM and INCLINE datasets and modify the right hand side vectors very slightly so that an exact solution is known for the purpose of the proceeding analysis. Let $\widehat{x}$ denote the approximate solution of the original problem $D^{-1}Ax=D^{-1}b$ to machine precision. If we now generate a slightly different right hand side vector $D^{-1}\widehat{b}=D^{-1}A\widehat{x}$, then $\widehat{x}$ becomes an exact solution for the 'nearby' rescaled linear system $D^{-1}A\widehat{x}=D^{-1}\widehat{b}$. For ease of notation, we will simply denote this linear system as $D^{-1}Ax=D^{-1}b$ and use $\widehat{x}$ to instead denote its approximate solution.

For a given stopping tolerance $\tau$, Tables \ref{tab:ferr} and \ref{tab:ferrb} show the exact relative forward error norm 
$\|\widehat{x} - x\|/\|x\|$, the estimated upper bound on the 
relative forward error norm ${\tt Ferr}$ and the normalized stopping tolerance $\epsilon$, both with and without row scaling, for the C2VSIM and INCLINE examples \footnote{${\tt Ferr}$ is computed from \eqref{eq:res} using the 
SuperLU \citep{SuperLU} routines {\tt dgscon}  
to estimate the condition number $\kappa(M^{-1}D^{-1}A^{-1})$.
}. 
With scaling, the upper bound on the relative 
forward error norm estimate is the same order of magnitude 
as $\epsilon$, although it is about an $\mathcal{O}(10^2)$ higher 
than the exact relative forward error norm. Without scaling, the 
estimate of the upper bound on the relative forward error norm is 
considerably larger than the normalized stopping tolerance $\epsilon$ and 
at least an $\mathcal{O}(10^5)$ larger than the exact relative forward error norm.

\begin{table}
\begin{center}
\begin{tabular}{|l||l|l|l||l|l|l|} \hline
 &  \multicolumn{3}{c||}{With row scaling} &  
 \multicolumn{3}{c|}{Without row scaling} \\ \hline
$log~\tau $  &  $\|\widehat{x} - x\|/\|x\|$  & ${\tt Ferr}$  & $\epsilon$
 & $\|\widehat{x} - x\|/\|x\|$  & ${\tt Ferr}$ & $\epsilon$\\
\hline
-1& 1.99E-4 & 6.24E-2 & 1.12E-2  & 2.96E-4 & 7.37E8 & 3.47E-11 \\
-2& 1.46E-5 & 1.81E-3 & 1.12E-3 & 1.76E-5 & 2.59E8 & 3.47E-12 \\
-3& 8.61E-7 & 2.30E-4 & 1.12E-4 & 9.49E-7 & 3.29E7 & 3.47E-13 \\
-4 & 5.09E-8 & 1.11E-5 &1.12E-5 & 1.27E-7 & 1.11E6 &  3.47E-14 \\
-5 & 5.09E-8 & 1.11E-5 &1.12E-6 & 1.05E-8 & 7.13E4 &  3.47E-15 \\
-6 & 1.48E-9 & 4.24E-7 &1.12E -7 & 1.05E-8 & 7.13E4 &  3.47E-16 \\
-7 & 1.09E-10 & 1.82E-8 &1.12E -8 & 1.06E-9 & 1.79E3 &  3.47E-17\\
-8 & 4.77E-12 & 9.51E-10 &1.12E -9 & 7.11E-11 & 2.03E2 &  3.47E-18\\
\hline
\end{tabular}
\caption{For a given stopping tolerance $\tau$, this Table compares 
the exact forward error norm, the estimated upper bound on the relative 
forward error norm and the normalized stopping tolerance $\epsilon$ from separately solving each of the linear 
systems $M^{-1}D^{-1}Ax=M^{-1}D^{-1}b$ (with row scaling) and 
$M^{-1}Ax=M^{-1}b$ (without row scaling) 
using PGMRES applied to the C2VSIM dataset.}  \label{tab:ferr}
\end{center}
\end{table}

\begin{table}
\begin{center}
\begin{tabular}{|l||l|l|l||l|l|l|} \hline
 &  \multicolumn{3}{c||}{With row scaling} &  
 \multicolumn{3}{c|}{Without row scaling} \\ \hline
$log~\tau $  &  $\|\widehat{x} - x\|/\|x\|$  & ${\tt Ferr}$  & $\epsilon$
 & $\|\widehat{x} - x\|/\|x\|$  & ${\tt Ferr}$ & $\epsilon$\\
\hline
-1& 2.56E-7 & 4.04E-5 &  3.83E-5 & 5.40E-7 & 1.11E-2 & 2.27E-8 \\
-2& 5.17E-8 & 3.55E-6 & 3.83E-6 & 3.87E-8 & 6.32E-3 & 2.27E-9 \\
-3& 4.01E-9 & 6.61E-7 & 3.83E-7 & 7.86E-9 & 3.86E-4 & 2.27E-10 \\
-4 & 4.61E-10 & 3.54E-8 & 3.83E-8 & 3.33E-10 & 1.07E-5 &  2.27E-11 \\
-5 & 3.49E-11 & 3.68E-9 & 3.83E-9 & 5.35E-11 & 1.05E-6 &  2.27E-12 \\
-6 & 6.97E-12 & 6.50E-10 & 3.83E-10 & 5.05E-12 & 7.37E-7 &  2.27E-13 \\
-7 & 2.27E-12 & 8.18E-11 &  3.83E-11 & 1.56E-12 & 1.37E-8 &  2.27E-14\\
-8 & 2.21E-12 & 3.80E-11 &  3.83E-12 & 1.50E-12 & 9.60E-9 &  2.27E-15\\
\hline
\end{tabular}
\caption{For a given stopping tolerance $\tau$, this Table compares 
the exact forward error norm, estimated upper bound on the relative 
forward error norm and the normalized stopping tolerance $\epsilon$ from separately solving each of the linear 
systems $M^{-1}D^{-1}Ax=M^{-1}D^{-1}b$ (with row scaling) and 
$M^{-1}Ax=M^{-1}b$ (without row scaling)
using PGMRES applied to the INCLINE dataset.}  \label{tab:ferrb}
\end{center}
\end{table}

\paragraph{Discussion} The main implications and limitations of the proceeding analysis are

\begin{itemize}
\item The coefficient matrices arising from integrated groundwater models should always be rescaled, irrespective of whether the stream nodes are represented in the coefficient matrix. This observation rests on the theory of forward error estimation, which is solver independent; 
\item The stopping tolerance in the GMRES algorithm should be chosen using the relation $\tau=\epsilon ||b||_*$, where $\epsilon$ is the desired target accuracy on the relative forward error. Even with $\epsilon$ fixed throughout a simulation, $\tau$ becomes a dynamical tolerance as $||b||_*$ varies in size; and
\item The empirically observed close correspondence between $\epsilon$ and {\tt Ferr} follows from the effectiveness of the preconditioner to reduce $\kappa(M^{-1}D^{-1}A)$ close to unity. The methodology described above is sufficiently broad that readers can follow our approach and establish correspondence for alternative choices of preconditioners and solvers. It also leads to the idea that a proxy can be found by comparing it with the exact error of the nearby rescaled linear system, without the need to estimate the forward error bound (to observe the effects of rescaling).
\end{itemize}



\section*{Implementation and Performance Benchmarking} \label{sect:perf}

Our Fortran 90 implementation of PGMRES(m) is adapted from the publically 
available sparse matrix package SPARSKIT \citep{Saad00}. This package 
provides general purpose Fortran 77 routines for preconditioning and 
iteratively solving sparse linear systems using matrices stored in 
compressed sparse row (CSR) format. We developed routines to efficiently 
convert IWFM storage format arrays into CSR format arrays and simultaneously rescale the linear system.

The PGMRES implementation uses a \emph{reverse communication} interface- an implementation style which enables the computationally intensive modules to be externally referenced. The sparse matrix-vector multiplication module and the LU solver are called at each iteration of PGMRES to form the Krylov subspace. In line 3 of the PGMRES algorithm, the vector $y:=Mx=AV_j$ is formed by sparse matrix-vector multiplication of $A$ and $v_j$. The LU solver then solves for a temporary vector $z$ using $y=\widehat{L}z$ followed by solving for $z$ using $z=\widehat{U}x$, where $\widehat{L}$ and $\widehat{U}$ are computed by the ILUT preconditioner prior to execution of the PGMRES algorithm. 

There are two distinct advantages to implementing an iterative solver with a reverse communication interface. Firstly, the PGMRES solver is independent of the choice of preconditioner and sparse matrix storage format because it does not directly process the matrices $A$ and $M$. Secondly, the bottleneck modules in the solver can be executed in high performance matrix algebra libraries which are tuned for the computer architecture. The matrix-vector multiply is classified as a sparse Basic Linear Algebra Subroutine (BLAS) level 2 operation and is implemented for a range of parallel computing and accelerator platforms including the newly emerging many-core CPU architectures. 

Our numerical experiments are performed using a Linux based Intel 
Fortran compiler V11.0 on a 2.00GHz Intel(R) Core(TM) 2 Duo CPU  (T6400) 
with 2MB cache. The relaxation parameter for the SOR method is set 
to $\omega=1.1$, the restart threshold of PGMRES is $m=20$ and 
the ILUT (ILU with threshold \citep{Saad00}) preconditioner has 
a drop tolerance\footnote{An element is replaced by zero if it is less 
than the drop tolerance multiplied by the original norm of the row 
containing the element.} of $0.01$ and maximum fill-in\footnote{Only the 
p largest elements in each upper and lower factor matrix are retained, 
the remainder are replaced by zero.} of $p=10$. We find by numerical 
experiment that this choice of PGMRES parameters gives optimal convergence 
rates for all of the non-symmetric datasets described in Table~\ref{tab:param}. 
In contrast, the optimal choice of $\omega$ varies significantly between each dataset. For the IWFM datasets, $1.1\leq\omega\leq 1.3$ is found to be an optimal range, whereas $\omega=1.95$ and $\omega=1.5$ are found to give optimal convergence rates for the NAC and INCLINE datasets. The sensitivity of the SOR method to $\omega$ is one reason to avoid using it for integrated groundwater flow models because of its sensitivity to changes in the linear systems.

The number of iterations and elapsed wall clock times for 
SOR and PGMRES($m$) applied to four of the datasets are respectively shown in 
Tables \ref{tab:bench:iter} and \ref{tab:bench:time} 
for a range of stopping tolerances $\tau$. The overall speedup from using PGMRES in place of the SOR solver is highly dependent on the dataset. For example, the speedup when using a tolerance of $1\times 10^{-8}$ is approximately 1.5x for the NAC dataset, 30x for the C2VSIM  and 240x for the BUTTE and INCLINE datasets. The corresponding approximation reduction in the number of iterations is respectively 6x, 90x, 800x and 900x. This reduction is attributed to the comparative effectiveness of the ILUT preconditioner at reducing the condition number of the preconditioned coefficient matrix and is a key feature necessary for a dynamical stopping tolerance $\tau=\epsilon ||D^{-1}b||$. In a separate experiment (not shown here), the GMRES method only needed about twice as many iterations to converge to machine precision. In contrast, the SOR method struggles to reach higher precision. The discrepancy between the realized speedups versus the reduction in the number of GMRES iterations is attributed to the cost of constructing the preconditioner. Once the preconditioner has been constructed, the relative additional cost of converging to higher precision is observed to correspond to the increase in the number of iterations.

\begin{table}
\begin{center}
 \begin{tabular}{|c||c|c||c|c||c|c||c|c|} \hline
  & \multicolumn{2}{c||}{C2VSIM} & 
    \multicolumn{2}{c||}{BUTTE} & 
    \multicolumn{2}{c||}{NAC} &
    \multicolumn{2}{c|}{INCLINE}  \\
\hline    
  $\log(\tau)$  &  SOR  & PGMRES & SOR  & PGMRES & SOR  & PGMRES & SOR & PGMRES\\
\hline
-1 &           41   & 5 & 65   & 6 & 166   & 46 &  167   & 26 \\
-2 &          155   & 6 &  310   & 7 & 219   & 54 &    785   & 30 \\
-3 &          276   & 6 &  2064   & 9 &  281   & 62 & 2918   & 33 \\
-4 &          397   & 7 & 3972   & 10 & 349   & 75 &  6275   & 37 \\
-5 &          518   & 8 & 5880   & 12 & 418   & 82 &  14523   & 41 \\
-6 &          638   & 9 & 7788   & 13 & 487   & 91 &  25835   & 46 \\
-7 &          759   & 9 & 9696   & 14 & 556   & 101 &  37148   & 49 \\
-8 &          880   & 10 & 11604   & 15 & 624   & 113 &  48461   & 52 \\
\hline
\end{tabular}
\caption{A comparison of the number of iterations
of SOR and PGMRES(m) applied to four of the datasets as the stopping 
tolerance $\tau$ is decremented.}
\label{tab:bench:iter}
\end{center}
\end{table}

\begin{table}
\begin{center}
 \begin{tabular}{|c||c|c||c|c||c|c||c|c|} \hline
  & \multicolumn{2}{c||}{C2VSIM} & 
    \multicolumn{2}{c||}{BUTTE} & 
    \multicolumn{2}{c||}{NAC} &
    \multicolumn{2}{c|}{INCLINE}  \\
 \hline
  $log(\tau)$  &  SOR  & PGMRES & SOR  & PGMRES & SOR  & PGMRES & SOR  & PGMRES \\
  \hline
-1 &  0.0244   &  0.00782  & 0.110        &  0.0478   &  0.601        & 0.696 & 0.36 & 0.253 \\
-2 &  0.0490   & 0.00787   & 0.510        &  0.0574   &  0.792        &  0.764 & 1.68 & 0.281 \\
-3 &   0.0883   & 0.00756   &  3.38        & 0.0572   &   1.02        & 0.868 & 6.22 & 0.293\\
-4 &    0.125        &  0.00801   &   6.47        & 0.0722   &  1.26        & 1.00 & 13.3 & 0.322\\
-5 &    0.162        &   0.00856   &  9.62        & 0.0661   & 1.50        & 1.11 & 31.0 & 0.348\\
-6 &   0.200        &  0.00898   & 12.7       &  0.0689   &  1.75       &  1.21 & 55.3 & 0.380\\
-7 &  0.238        & 0.00922   &  15.8        & 0.0769   &   2.00        & 1.33 & 79.2 & 0.407\\
-8 &     0.276        & 0.00963   &  18.9        &  0.0792   & 2.24        &  1.48 & 103.0& 0.434\\
\hline
 \end{tabular}
\caption{A comparison of elapsed wall clock time (in seconds) 
of SOR and PGMRES(m) applied to four of the datasets as the stopping 
tolerance $\tau$ is decremented. 
All timings are reported to three significant figures.}
\label{tab:bench:time}
\end{center}
 \end{table}

Finally, Table 6 shows the overall performance 
improvement in the IWFM simulation using the PGMRES(m) solver 
in place of the SOR solver and the proportion (shown in parentheses) of overall computation 
spent in the preconditioner and solvers for three of the datasets. 
The C2VSIM and C2VSIM9 simulations are run over 82 years at monthly 
increments (984 time steps) and the HCMP simulation is run over 2 years 
at weekly increments (104 time steps). The normalized tolerance $\epsilon$ in the linear solver was fixed throughout the simulation at $1\times 10^{-5}$ which is $0.1$ times the relative data accuracy of $1\times 10^{-4}$.

As expected, the overall performance gains from using PGMRES(m) are more prominent with the larger datasets since the absolute time reduction is most significant- HCMP and C2VSIM9 exhibit 7.74x and 7.56x speedups respectively. The reduction of the total time spent in the solver is also considerably more significant for the C2VSIM9 dataset than the HCMP and C2VSIM datasets - 64.3x reduction versus  14.4x and 19.0x respectively.
This relative improvement is different from the isolated benchmarking results (some of which are shown in Tables \ref{tab:bench:iter} and \ref{tab:bench:time},) in which the speedup from using PGMRES is 150x, 10x and 30x when using the C2VSIM9, HCMP and C2VSIM datasets respectively. This indicates that the linear system used for the linear solver benchmarking in Tables \ref{tab:bench:iter} and \ref{tab:bench:time} leads to optimistic performance gains compared to the total linear solver time reduction over the simulation, in which the linear systems may change significantly. With PGMRES, iterative solution of the linear system ceases to be a major bottle neck for the C2VSIM and C2VSIM9 datasets - only 9.12\% and 9.65\% of the overall simulation time is spent in the linear solver. This is the reason why the IWFM speedups of 2.2x and 7.56x are significantly lower than the total solver time reductions of 19.0x and 64.3x. Conversely, the linear solver is still a major bottleneck in IWFM with the HCMP dataset accounting for 45.1\% of the IWFM simulation time.  
 
\begin{table}
\begin{center}
\begin{tabular}{|c|c|c|c|}
\hline
& HCMP & C2VSIM & C2VSIM9\\
\hline
IWFM (SOR) & 20.4 (84.0\%) & 15.59 (79.0\%) & 121 (82.0\%) \\
IWFM (PGMRES)& 2.63 (45.1\%) & 7.12 (9.12\%) & 16.0 (9.65\%) \\
\hline
Speedup& 7.74x (14.4x) & 2.2x (19.0x) & 7.56x (64.3x)  \\
\hline
\end{tabular}
\caption{This Table shows the time in minutes and proportion of IWFM simulation time (in parenthesis) spent in the solvers for each of the datasets. The bottom row shows the speedup in IWFM simulation time and total solver time (in parenthesis) if the SOR solver is replaced by PGMRES. The C2VSIM and C2VSIM9 simulations use 984 time steps and the HCMP simulation uses 104 time steps. All timings are reported to three significant figures.}
\label{tab:bench_final}
\end{center}
\end{table}

\section*{Conclusion} \label{sect:conclusion}

An open problem that arises in many modern iterative linear solvers is how to choose the residual tolerance in the linear solver to be consistent with the tolerance on the solution error. This article firstly illustrates the scaling issues arising in linear systems from integrated groundwater and surface water models. It is shown that re-scaling is a key step towards control of the 
forward error, irrespective of the choice of solver and preconditioner. We implemented a preconditioned GMRES algorithm and observe that the normalized stopping tolerance on the estimated residual norm becomes a proxy for the estimated upper bound on the forward error when the linear systems are rescaled. Furthermore, the rescaling and modified error control are simple to implement at negligible additional computational cost.

This article demonstrates a number of favorable properties of PGMRES: (i) the optimal ILUT preconditioner parameter settings are independent of the dataset; (ii) it is well suited to adaptive residual error control; (iii) when benchmarked against the SOR method, the comparative speedups with rescaling and error control can lead to significant overall application speedups (as high as 7.7x for IWFM); and (iv)  performance profiling shows that the new linear solver removes a major performance bottleneck in IWFM.

This research is expected to broadly impact groundwater modelers by demonstrating a practical and general approach for setting the residual tolerance in line with the solution error tolerance.

\nolinenumbers 
\section*{Acknowledgments}
The authors would like to thank Dr Richard Niswonger for his support in investigating the use of Newton-GMRES methods for MODFLOW-2005 and providing the INCLINE dataset. Any inaccuracies in statements concerning MODFLOW-2005 or related USGS software packages are solely the responsibility of the authors. The authors also thank the Texas Water Development Board for providing the report and model files for the NAC dataset.


\end{document}